\documentclass{ifacconf}
\usepackage{graphicx}      
\usepackage{natbib}        

\usepackage{graphics} 

\usepackage{times} 

\usepackage{tikz} 
\usepackage{pgfplots} 

\usepackage{algpseudocode}%
\usepackage{algorithm}
\usepackage{amsmath}
\usepackage{amssymb}
\usepackage{epstopdf}

\usepackage{color}
\def\tr{\text{tr}}

\usepackage{lipsum}
\usepackage{verbatim}

\usepackage{mathtools}
\usepackage{times}

\newtheorem{theorem}{\bf Theorem}  
 \newtheorem{remark}[theorem]{\bf Remark}
 \newtheorem{proposition}[theorem]{\bf Proposition} 
\newtheorem{assumption}[theorem]{\bf Assumption}

\begin{document}
\begin{frontmatter}

\title{Data-driven analysis and control of continuous-time systems under aperiodic sampling\thanksref{footnoteinfo}} 

\thanks[footnoteinfo]{This work was funded by Deutsche Forschungsgemeinschaft (DFG, German Research Foundation) under Germany’s Excellence Strategy - EXC 2075 - 390740016 and under grant AL 316/13-2 - 285825138.
The authors thank the International Max Planck Research School for Intelligent Systems (IMPRS-IS) for supporting Julian Berberich.\\
\textcopyright 2021 the authors. This work has been accepted to IFAC for publication under a 
Creative Commons Licence CC-BY-NC-ND.}

\author[First]{Julian Berberich} 
\author[First]{Stefan Wildhagen} 
\author[First]{Michael Hertneck}
\author[First]{Frank Allg\"ower}

\address[First]{University of Stuttgart, Institute for Systems Theory and Automatic Control, 70569 Stuttgart, Germany (e-mail: \{julian.berberich, stefan.wildhagen, michael.hertneck, frank.allgower\}@ist.uni-stuttgart.de).}

\begin{abstract}                
We investigate stability analysis and controller design of unknown continuous-time systems under state-feedback with aperiodic sampling, using only noisy data but no model knowledge.
We first derive a novel data-dependent parametrization of all linear time-invariant continuous-time systems which are consistent with the measured data and the assumed noise bound.
Based on this parametrization and by combining tools from robust control theory and the time-delay approach to sampled-data control, we compute lower bounds on the maximum sampling interval (MSI) for closed-loop stability under a given state-feedback gain, and beyond that, we design controllers which exhibit a possibly large MSI.
Our methods guarantee the stability properties robustly for all systems consistent with the measured data.
As a technical contribution, the proposed approach embeds existing methods for sampled-data control into a general robust control framework, which can be directly extended to model-based robust controller design for uncertain time-delay systems under general uncertainty descriptions.
\end{abstract}


\end{frontmatter}

\section{Introduction}
For a given continuous-time plant $\dot{x}(t)=Ax(t)+Bu(t)$ and a sampled-data state-feedback controller $u(t)=Kx(t_k),t\in[t_k,t_{k+1})$, we define the maximum sampling interval (MSI) as the largest possible $h>0$ such that the closed loop is asymptotically stable for any sampling sequence with $t_{k+1}-t_k\leq h$.
Computing possibly large lower bounds on the MSI is an interesting and relevant research question since its knowledge allows to reduce the sampling frequency as much as possible without affecting stability.
Furthermore, the stability analysis of sampled-data systems can be applied to many practical problems, e.g., to study the behavior of networked control systems with transmission delays or packet dropouts.
Since the MSI is an essential concept in sampled-data control (compare~\cite{hetel2017recent,hespanha2007survey}), a wide variety of methods to compute such MSI bounds have been developed in the literature, using different viewpoints such as the time-delay approach (\cite{fridman2004robust,fridman2010refined}), the hybrid system approach (\cite{carnevale2007lyapunov,nesic2009explicit}), the discrete-time approach (\cite{fujioka2009stability,hetel2011delay}), and the input/output stability approach (\cite{mirkin2007some}), which are all discussed extensively in the survey~\cite{hetel2017recent}.
Typically, these methods are model-based, i.e., they require an (accurate) model for their implementation.
While measured data of a system is typically simple to obtain, deriving a model via first principles can be complicated and time-consuming.
If no model is available, \emph{finding} a controller with guaranteed stability for all sampling intervals in a desired range is particularly important since, without model knowledge, it is not even clear how to design a stabilizing controller under continuous state-feedback.
Nevertheless, to the best of our knowledge, computing guaranteed MSI bounds and designing sampled-data controllers based only on measured data has not been considered in the literature thus far.

In this paper, we develop a direct approach for a) computing MSI bounds for unknown systems with a given controller and b) designing controllers with possibly large MSI bounds, both based only on a noisy data trajectory of finite length.
Our approach is inspired by recent results of~\cite{waarde2020from} on discrete-time data-driven control, which are based on a data-dependent parametrization of all matrices consistent with the measured data and a given noise bound.
Our work is also closely related to~\cite{berberich2020combining}, which uses similar ideas for combining data with additional prior knowledge for robust controller design.
In the present paper, we first derive a result similar to~\cite{waarde2020from} to parametrize unknown \emph{continuous-time} systems using one a priori collected trajectory.
We combine this parametrization with well-established robust control techniques (\cite{scherer2000linear,scherer2001lpv}) to verify existing sampled-data stabilization conditions of~\cite{fridman2010refined} robustly for all systems which are consistent with the measured data and the assumed noise bound.
While most existing results for sampled-data control of uncertain systems using time-delay arguments only apply to restricted classes of uncertainties, our results lead to a general robust sampled-data control formulation with larger flexibility.
Another interpretation of our work is that of providing tools for analyzing unknown systems w.r.t. a particular sampled-data system property (the MSI) using only measured data, similar to other recent works on data-driven system analysis such as, e.g.,~\cite{maupong2017lyapunov,koch2020provably}.

A well-established alternative to the results presented in this paper is to use the measured data to identify the unknown system and then apply model-based tools as presented in~\cite{hetel2017recent} to the estimated system.
Our approach is more direct and does not rely on an intermediate estimation step while at the same time being non-conservative w.r.t. the model-based conditions in~\cite{fridman2010refined}.
On the other hand, providing tight estimation bounds using noisy data of finite length is in general non-trivial for stochastic noise (\cite{matni2019tutorial,matni2019self}) and can be computationally expensive for deterministic noise (\cite{milanese1991optimal}).
Another closely related paper by~\cite{rueda2020delay} on data-driven control for time-delay systems using delayed measurements addresses a conceptually similar control goal (stability under uncertain delays) but considers a discrete-time scenario which leads to an entirely different problem formulation.

\subsection*{Notation}
For a symmetric matrix $P$, we write $P\succ0$ ($P\succeq0$) or $P\prec0$ ($P\preceq0$) if $P$ is positive (semi-) definite or negative (semi-) definite, respectively.
Further, in matrix inequalities, $\star$ denotes entries which can be inferred from symmetry.
We denote the set of nonnegative integers by $\mathbb{N}_0$.

\section{Preliminaries}
We consider an unknown linear time-invariant (LTI) system
\begin{align}\label{eq:sys_cont}
\dot{x}(t)=A_{\tr}x(t)+B_{\tr}u(t),\>\>x(0)=x_0,
\end{align}
with $x_0,x(t)\in\mathbb{R}^n$, $u(t)\in\mathbb{R}^m$, $t\geq0$, where the true system matrices $A_{\tr},B_{\tr}$ are unknown.
Our goal is to design a sampled-data state-feedback controller $u(t)=Kx(t_k)$ for $t\in[t_k,t_{k+1})$, where $\{t_k\}_{k=0}^\infty$ with $t_0=0$ and $\lim_{k\to\infty}t_k=\infty$ are the sampling times, such that the closed loop is exponentially stable for all sampling times satisfying $t_{k+1}-t_k\in(0,h],k\in\mathbb{N}_0$, with a possibly large $h>0$.
We call the maximum $h$ such that the sampled-data system is exponentially stable the MSI.
We assume that $N\in\mathbb{N}_0$ measurements $\{\dot{x}(\tau_k),x(\tau_k),u(\tau_k)\}_{k=1}^N$, sampled at not necessarily equidistant sampling times $\{\tau_k\}_{k=1}^N$, of the perturbed system
\begin{align}\label{eq:sys_cont_dist}
\dot{x}(t)=A_{\tr}x(t)+B_{\tr}u(t)+B_dd(t),\>\>x(0)=x_0,
\end{align}
for some unknown disturbance signal $d(t)\in\mathbb{R}^{m_d}$, $t\geq0$, are available.
The matrix $B_d\in\mathbb{R}^{n\times m_d}$ has full column rank and is known which allows us to model the influence of the disturbance on the system.
We define the data matrices
\begin{align*}
\dot{X}&\coloneqq\begin{bmatrix}\dot{x}(\tau_1)&\dot{x}(\tau_2)&\dots&\dot{x}(\tau_N)\end{bmatrix},\\
X&\coloneqq\begin{bmatrix}x(\tau_1)&x(\tau_2)&\dots&x(\tau_N)\end{bmatrix},\\
U&\coloneqq\begin{bmatrix}u(\tau_1)&u(\tau_2)&\dots&u(\tau_N)\end{bmatrix}.
\end{align*}
We assume that the (unknown) disturbance sequence affecting the measured data, denoted by $\{\hat{d}(\tau_k)\}_{k=1}^N$, satisfies a bound of the form $\hat{D}\coloneqq\begin{bmatrix}\hat{d}(\tau_1)&\hat{d}(\tau_2)&\dots&\hat{d}(\tau_N)\end{bmatrix}\in\mathcal{D}$ with
\begin{align}\label{eq:noise_bound}
\mathcal{D}\coloneqq\Big\{D\in\mathbb{R}^{m_d\times N}\Bigm|
\begin{bmatrix}D^\top\\I\end{bmatrix}^\top
\begin{bmatrix}Q_d&S_d\\S_d^\top&R_d\end{bmatrix}
\begin{bmatrix}D^\top\\I\end{bmatrix}\succeq0\Big\}
\end{align}
for some $Q_d\in\mathbb{R}^{N\times N}$, $S_d\in\mathbb{R}^{N\times m_d}$, $R_d\in\mathbb{R}^{m_d\times m_d}$ with $Q_d\prec0$.
This includes a wide range of practically relevant noise bounds (compare~\cite{persis2020formulas,berberich2020design,waarde2020from,koch2020provably,berberich2020combining} where similar noise bounds are considered in a discrete-time setting).
A key practical issue in obtaining data as above is that time derivatives $\dot{X}$ of the state are typically not accessible directly.
However, after having recorded a densely-sampled trajectory of the state $x$, one could obtain a good approximation of $\dot{x}$ in a post-processing step, e.g., using finite differences or even a non-causal filter.
Alternatively, various methods for derivative estimation have been developed in the field of continuous-time system identification, compare, e.g.,~\cite{garnier2003continuous}.
Since the data are recorded in an offline experiment, the mentioned approaches are also amenable for networked control systems:
Therein, a sensor located at the plant could measure the state with a high sampling rate and store the data, which would then be extracted after the experiment. 
Note that the network is only required for closed-loop operation, such that the achievable sampling rate of the experiment depends only on the sensor and not on the network.

Further, we do not assume that $\dot{x}(\tau_k)$ can be measured exactly.
In particular, including the disturbance $d$ in~\eqref{eq:sys_cont_dist} allows for measurement noise on $\dot{x}(\tau_k)$.
For instance, if we assume that the deviation of the measured value of $\dot{x}(\tau_k)$ from the true value is bounded in the $2$-norm by $\bar{d}>0$ for each $k=1,\dots,N$, then we can choose $B_d=I$, $Q_d=-I$, $S_d=0$, $R_d=\bar{d}^2NI$.
In Appendix A, we show that an estimate of $\dot{x}(\tau_k)$ as well as an error bound can be obtained from measured data via Euler discretization if norm bounds on the matrices $A_{\tr}$ and $B_{\tr}$ are available.
It is an interesting issue for future research to relax the latter assumption, extend this approach to more general discretization schemes, and validate it with a realistic example.

Let us now define the set of all systems consistent with the measured data and the noise bound
\begin{align*}
\Sigma_c\coloneqq\{(A,B)\mid \dot{X}=AX+BU+B_dD,D\in\mathcal{D}\}.
\end{align*}
On a technical level, our contribution can be summarized as follows:
1) We derive a simple, data-dependent parametrization of the set $\Sigma_c$ (Section~\ref{sec:param}), which we then employ to 2) \emph{compute} a lower bound on the MSI under a given state-feedback (Section~\ref{sec:analysis}) and to 3) \emph{design} a state-feedback controller for~\eqref{eq:sys_cont} with a possibly large MSI (Section~\ref{sec:design}).
The MSI bounds are guaranteed robustly for all system matrices $\begin{bmatrix}A&B\end{bmatrix}\in\Sigma_c$ which are consistent with the measured data and the noise bound.
Thus, our results provide an alternative to model-based control under aperiodic sampling as surveyed in~\cite{hetel2017recent}, requiring no model knowledge but only measured data.

Finally, we note that it is possible to extend the results of this paper to output-feedback control under aperiodic sampling by considering an extended state vector based on the first $n$ time derivatives of the input and output, similar to the discrete-time output-feedback formulation in~\cite{berberich2020combining}.

\section{Data-driven system parametrization}\label{sec:param}
In this section, we provide a purely data-driven parametrization of all continuous-time LTI systems as in~\eqref{eq:sys_cont} which are consistent with the measured data and the noise bound.
The presented approach provides a variation of~\cite{waarde2020from}, where analogous results are obtained for discrete-time systems.
To this end, we define $Z\coloneqq\begin{bmatrix}X\\U\end{bmatrix}$ as well as
\begin{align}\label{eq:bound_def}
P_c=\begin{bmatrix}Q_c&S_c\\S_c^\top&R_c\end{bmatrix}\coloneqq\begin{bmatrix}-Z&0\\\dot{X}&B_d\end{bmatrix}
\begin{bmatrix}Q_d&S_d\\S_d^\top&R_d\end{bmatrix}
\begin{bmatrix}-Z&0\\\dot{X}&B_d\end{bmatrix}^\top.
\end{align}
\begin{theorem}\label{thm:cont_param}
It holds that
\begin{align}\label{eq:thm_cont_param}
\Sigma_c=\Big\{(A,B)\Bigm|
\begin{bmatrix}\begin{bmatrix}A&B\end{bmatrix}^\top\\I\end{bmatrix}^\top
P_c
\begin{bmatrix}\begin{bmatrix}A&B\end{bmatrix}^\top\\I\end{bmatrix}\succeq0\Big\}.
\end{align}
\end{theorem}
\textbf{Proof}$\quad$
The proof follows similar arguments as~\cite{waarde2020from} and hence we only give a sketch for the special case $B_d=I$ and refer to~\cite{waarde2020from} for details.
If $B_d=I$, then $(A,B)$ lies in the set on the right-hand side of~\eqref{eq:thm_cont_param} if and only if
\begin{align}\label{eq:thm_cont_param_proof1}
\begin{bmatrix}(\dot{X}-AX-BU)^\top\\I\end{bmatrix}^\top
\begin{bmatrix}Q_d&S_d\\S_d^\top&R_d\end{bmatrix}
\begin{bmatrix}\star\\\star\end{bmatrix}\succeq0.
\end{align}
It is straightforward to see that~\eqref{eq:thm_cont_param_proof1} is in turn equivalent to the existence of $D\in\mathcal{D}$ satisfying $\dot{X}=AX+BU+D$, thus concluding the proof.
$\hfill\blacksquare$

Theorem~\ref{thm:cont_param} provides a direct data-driven parametrization of all continuous-time systems consistent with the data using noisy input-state measurements, which need not be sampled at equidistant sampling instants.
The proof is analogous to that of the discrete-time results by~\cite{waarde2020from} and it relies on simple algebraic manipulations of the data equation.
The system matrices $(A,B)$ can be interpreted as an \emph{uncertainty} satisfying the quadratic matrix inequality bound~\eqref{eq:thm_cont_param}.
Indeed, the remainder of this paper makes extensive use of this interpretation by applying S-procedure-based linear matrix inequality (LMI) relaxations to robustify existing model-based LMIs for computing MSI bounds against the uncertainty class $\Sigma_c$.

In~\cite{persis2020formulas}, an alternative approach to data-driven control of continuous-time systems is outlined.
Controller design based on Theorem~\ref{thm:cont_param} has multiple advantages if compared to the approach by~\cite{persis2020formulas} such as guaranteed robustness with respect to noise and a lower computational complexity (the number of decision variables for controller design in the present paper is independent of the data length).
It is an interesting issue for future research to employ Theorem~\ref{thm:cont_param} in order to develop further continuous-time data-driven control methods beyond the sampled-data formulations considered in the present paper.

\section{Stability analysis under aperiodic sampling}\label{sec:analysis}
In this section, we analyze stability of unknown sampled-data systems and compute lower bounds on the MSI by modeling the closed-loop system as a time-delay system.
To this end, our work builds on the model-based time-delay approach to sampled-data control which is discussed, e.g., in~\cite{hetel2017recent,fridman2004robust,fridman2010refined,liu2012wirtingers}.
The key idea of the time-delay approach to sampled-data control is to write the controller as
\begin{align}\label{eq:sampled_data_controller}
u(t)=Kx(t_k)=Kx(t-\tau(t)),\quad
\tau(t)\coloneqq t-t_k
\end{align}
for all $t\in[t_k,t_{k+1})$, where $\tau$ is a delay with $\tau(t_k)=0$ and $\dot{\tau}(t)=1$ if $t\neq t_k$.
The closed-loop sampled-data control system can then be written as
\begin{align}\label{eq:sys_delay}
\dot{x}(t)&=Ax(t)+BKx(t-\tau(t))
\end{align}
for any $t\geq0$, i.e., as a linear time-delay system.
The literature contains various methods to compute values $h$ such that~\eqref{eq:sys_delay} is exponentially stable with time-varying delay $\tau(t)$ for any $\tau(t)\in(0,h]$, see e.g.~\cite{hetel2017recent,fridman2004robust,fridman2010refined}.
Throughout this paper, we focus on the conditions stated in~\cite[Proposition 1]{fridman2010refined} which rely on the Lyapunov-Krasovskii functional
\begin{align}\label{eq:LKF}
V(t,x(t),\dot{x}(t))=&\>x(t)^\top P_1x(t)\\\nonumber
&+(h-\tau(t))\int_{t-\tau(t)}^t \dot{x}(s)^\top R\dot{x}(s)\text{d}s,
\end{align}
where $P_1\succ0,R\succ0,\tau(t)=t-t_k$.
Using~\eqref{eq:LKF},~\cite{fridman2010refined} computes an MSI bound based on model knowledge as follows.
\begin{proposition}\label{prop:fridman}
(\cite{fridman2010refined})
The system $\dot{x}=Ax(t)+Bu(t)$ with controller~\eqref{eq:sampled_data_controller} is exponentially stable for any sampling sequence $\{t_k\}_{k=0}^{\infty}$ with $t_{k+1}-t_k\in(0,h]$, $t_0=0$, $\lim_{k\to\infty}t_k=\infty$, if there exist $R\succ0$, $P_1\succ0$, $P_2$ and $P_3$ such that
\begin{align}\label{eq:time_delay1}
&\begin{bmatrix}P_2^\top(A+BK)+(A+BK)^\top P_2&\star\\
P_1-P_2+P_3^\top (A+BK)&-P_3-P_3^\top+hR
\end{bmatrix}\prec0,\\\label{eq:time_delay2}
&\begin{bmatrix}P_2^\top(A+BK)+(A+BK)^\top P_2&\star&\star\\
P_1-P_2+P_3^\top(A+BK)&-P_3-P_3^\top&\star\\
-h(BK)^\top P_2&-h(BK)^\top P_3&-hR
\end{bmatrix}\prec0.
\end{align}
\end{proposition}
Throughout this paper, we make the following technical assumption on the matrix $P_c$ involved in the data-dependent parametrization provided by Theorem~\ref{thm:cont_param}.
\begin{assumption}\label{ass:inertia}
The matrix $P_c$ is invertible and has exactly $m_d$ positive eigenvalues.
\end{assumption}
This assumption is not restrictive and it holds if the available measurements are sufficiently rich, i.e., the matrix $Z$ has full row rank, under mild additional assumptions.
To be precise, we show in Appendix B that Assumption~\ref{ass:inertia} holds for the common special case $S_d=0$ if (i) $Z$ has full row rank, (ii) $B_d$ is invertible, i.e., $m_d=n$, and (iii) the disturbance $\hat{D}$ generating the data satisfies a strict version of the quadratic matrix inequality in~\eqref{eq:noise_bound}, i.e., $\hat{D}Q_d\hat{D}^\top+R_d\succ0$.
Note that $Z$ having full row rank does not imply that the underlying system can be uniquely identified from the available data due to the influence of the disturbance.
Defining 
\begin{align*}
\begin{bmatrix}\tilde{Q}_c&\tilde{S}_c\\\tilde{S}_c^\top&\tilde{R}_c\end{bmatrix}\coloneqq
\begin{bmatrix}Q_c&S_c\\S_c^\top&R_c\end{bmatrix}^{-1},\>
\tilde{P}_c\coloneqq\begin{bmatrix}-\tilde{R}_c&\tilde{S}_c^\top\\\tilde{S}_c&-\tilde{Q}_c
\end{bmatrix},
\end{align*}
we now provide a data-based stability certificate for the considered sampled-data system.

\begin{theorem}\label{thm:analysis}
If Assumption~\ref{ass:inertia} holds, then $\dot{x}(t)=Ax(t)+Bu(t)$ with controller~\eqref{eq:sampled_data_controller} is exponentially stable for any $\begin{bmatrix}A&B\end{bmatrix}\in\Sigma_{c}$ and any sampling sequence $\{t_k\}_{k=0}^{\infty}$ with $t_{k+1}-t_k\in(0,h]$, $t_0=0$, and $\lim_{k\to\infty}t_k=\infty$, if there exist ${\lambda_1>0}$, ${\lambda_2>0}$, ${P_1\succ0}$, $R\succ0$, $P_2$, and $P_3$ such that $P=\begin{bmatrix}P_1&0\\P_2&P_3\end{bmatrix}$, ${P_R=\text{diag}(P,R)}$, $P_{R2}=\begin{bmatrix}P_1&0\\P_2&P_3\\0&R\end{bmatrix}$ and
\begin{align}\label{eq:thm_time_delay_LMI1}
&\left[
\begin{array}{cc}
I&0\\
\begin{bmatrix}0&I\\0&-I\\0&\frac{1}{2}hI\end{bmatrix}&\begin{bmatrix}0\\I\\0\end{bmatrix}\\\hline
0&I\\
\begin{bmatrix}I&0\\K&0\end{bmatrix}&\begin{bmatrix}0\\0\end{bmatrix}
\end{array}
\right]^\top
\left[
\begin{array}{c|c}
\begin{matrix}
0&P_{R2}^\top\\P_{R2}&0
\end{matrix}&0\\\hline
0&\lambda_1\tilde{P}_c
\end{array}
\right]
\left[
\begin{array}{cc}
\star&\star\\\star&\star\\\hline
\star&\star\\\star&\star
\end{array}
\right]\prec0,
\end{align}
\begin{align}\label{eq:thm_time_delay_LMI2}
&\left[
\begin{array}{cc}
I&0\\
\begin{bmatrix}0&I&0\\0&-I&0\\0&0&-\frac{1}{2}hI\end{bmatrix}&\begin{bmatrix}0\\I\\0\end{bmatrix}\\\hline
0&I\\
\begin{bmatrix}I&0&0\\K&0&-hK\end{bmatrix}&\begin{bmatrix}0\\0\end{bmatrix}
\end{array}
\right]^\top
\left[
\begin{array}{c|c}
\begin{matrix}
0&P_{R}^\top\\P_{R}&0
\end{matrix}&0\\\hline
0&\lambda_2\tilde{P}_c
\end{array}
\right]
\left[
\begin{array}{cc}
\star&\star\\\star&\star\\\hline
\star&\star\\\star&\star
\end{array}
\right]\prec0.
\end{align}
\end{theorem}
\textbf{Proof}$\quad$
Applying the dualization lemma~\cite[Lemma 4.9]{scherer2000linear} to~\eqref{eq:thm_cont_param} (the required inertia properties hold by Assumption~\ref{ass:inertia}), we have that $\begin{bmatrix}A&B\end{bmatrix}\in\Sigma_c$ if and only if
\begin{align}\label{eq:AB_bound_primal}
\begin{bmatrix}\begin{bmatrix}A&B\end{bmatrix}\\I\end{bmatrix}^\top
\begin{bmatrix}-\tilde{R}_c&\tilde{S}_c^\top\\\tilde{S}_c&-\tilde{Q}_c\end{bmatrix}
\begin{bmatrix}\begin{bmatrix}A&B\end{bmatrix}\\I\end{bmatrix}\succeq0.
\end{align}
We note that this step is analogous to the equivalence of~\cite[Conditions (4.4.1) and (4.4.4)]{scherer2000linear}, where the strict and non-strict inequalities can be interchanged.
Note that
\begin{align*}
\begin{bmatrix}0&I\\A+BK&-I\\0&\frac{1}{2}h\end{bmatrix}=
\begin{bmatrix}0&I\\0&-I\\0&\frac{1}{2}hI\end{bmatrix}+\begin{bmatrix}0\\I\\0\end{bmatrix}
\begin{bmatrix}A&B\end{bmatrix}\begin{bmatrix}I&0\\K&0\end{bmatrix}.
\end{align*}
Therefore, using the full-block S-procedure~(\cite{scherer2001lpv}), \eqref{eq:thm_time_delay_LMI1} implies that for any $\begin{bmatrix}A&B\end{bmatrix}$ satisfying~\eqref{eq:AB_bound_primal} the matrix
\begin{align*}
\begin{bmatrix}P_1&0\\P_2&P_3\\0&R\end{bmatrix}^\top\begin{bmatrix}0&I\\A+BK&-I\\0&\frac{1}{2}hI\end{bmatrix}
+\begin{bmatrix}0&I\\A+BK&-I\\0&\frac{1}{2}hI\end{bmatrix}^\top \begin{bmatrix}P_1&0\\P_2&P_3\\0&R\end{bmatrix}
\end{align*}
is negative definite, i.e.,~\eqref{eq:time_delay1} holds.
Using that~\eqref{eq:AB_bound_primal} holds if and only if $\begin{bmatrix}A&B\end{bmatrix}\in\Sigma_c$, this in turn implies that~\eqref{eq:time_delay1} holds for any $\begin{bmatrix}A&B\end{bmatrix}\in\Sigma_c$.
Similarly, we can show that~\eqref{eq:thm_time_delay_LMI2} implies
\begin{align*}
&\begin{bmatrix}P_1&0&0\\P_2&P_3&0\\0&0&R\end{bmatrix}^\top
\begin{bmatrix}0&I&0\\A+BK&-I&-hBK\\0&0&-\frac{1}{2}hI\end{bmatrix}\\
+&\begin{bmatrix}0&I&0\\A+BK&-I&-hBK\\0&0&-\frac{1}{2}hI\end{bmatrix}^\top\begin{bmatrix}P_1&0&0\\P_2&P_3&0\\0&0&R\end{bmatrix}\prec0
\end{align*}
for any $\begin{bmatrix}A&B\end{bmatrix}\in\Sigma_c$, i.e.,~\eqref{eq:time_delay2} holds.
Thus, by Proposition~\ref{prop:fridman}, we can deduce exponential stability for any $\begin{bmatrix}A&B\end{bmatrix}\in\Sigma_c$.
$\hfill\blacksquare$

Theorem~\ref{thm:analysis} provides a purely data-dependent approach to verify closed-loop stability of a sampled-data system with aperiodic sampling.
For any fixed MSI bound $h$, the conditions~\eqref{eq:thm_time_delay_LMI1}--\eqref{eq:thm_time_delay_LMI2} are LMIs and thus, a possibly large $h$ leading to closed-loop robust stability can be found via a simple bisection algorithm.
The derived stability guarantees hold for an \emph{arbitrary} sampling sequence $\{t_k\}_{k=0}^{\infty}$ of infinite length satisfying $t_0=0$, $t_{k+1}-t_k\in(0,h]$, and $\lim_{k\to\infty}t_k=\infty$.
This only requires availability of a noisy open-loop data trajectory of finite length which is generated by a sampling sequence $\{\tau_k\}_{k=1}^N$ \emph{not} necessarily satisfying $\tau_{k+1}-\tau_k\leq h$.
Under mild technical assumptions, the statement is non-conservative in the sense that the proposed conditions are equivalent to~\eqref{eq:time_delay1}-\eqref{eq:time_delay2} holding for all $\begin{bmatrix}A&B\end{bmatrix}\in\Sigma_{c}$.
This is due to the fact that the employed S-procedure only uses one matrix inequality constraint~\eqref{eq:AB_bound_primal}, for which relaxations with a scalar multiplier are exact, compare~\cite{waarde2020from} for details.
The proof relies on an application of robust control arguments from~\cite{scherer2000linear,scherer2001lpv} to verify the model-based LMIs~\eqref{eq:time_delay1}-\eqref{eq:time_delay2} robustly for all systems consistent with the measured data.
Since the parametrization of $\Sigma_c$ in Theorem~\ref{thm:cont_param} takes a ``dual'' form (involving $\begin{bmatrix}A&B\end{bmatrix}^\top$ instead of $\begin{bmatrix}A&B\end{bmatrix}$), we use the dualization lemma to derive the ``primal'' parametrization~\eqref{eq:AB_bound_primal} which we then combine with Inequalities~\eqref{eq:time_delay1}-\eqref{eq:time_delay2}.
It is worth noting that, to robustify Inequality~\eqref{eq:time_delay1} w.r.t. $\begin{bmatrix}A&B\end{bmatrix}\in\Sigma_c$, we use as a technical argument a \emph{non-square} matrix $P_{R2}$ which plays the role of the Lyapunov matrix in robust control.
Finally, we conjecture that an extension of Theorem~\ref{thm:design} to include performance criteria for the closed loop is straightforward, e.g., by robustly verifying the performance conditions in~\cite{fridman2010refined}.

\section{Controller design for closed-loop stability under aperiodic sampling}\label{sec:design}
While Theorem~\ref{thm:analysis} provides an effective way to compute MSI bounds of an unknown plant, the result requires that a state-feedback gain $K$ is given.
In case no controller is available, a matrix $K$ with $A+BK$ Hurwitz for all $(A,B)\in\Sigma_c$ can be found by combining Theorem~\ref{thm:cont_param} with a Lyapunov LMI and the S-procedure, similar to Theorem~\ref{thm:analysis} or the discrete-time results by~\cite{waarde2020from}.
However, applying the analysis conditions in Theorem~\ref{thm:analysis} for the resulting controller may lead to a small MSI such that stability can only be achieved under an unnecessarily high sampling rate.
In particular, if achieving a small sampling rate constitutes the main control objective, then we want to find a possibly large MSI bound $h$ holding for \emph{some} controller, which needs to be found concurrently with $h$.
Since the matrix $\tilde{Q}_c$ is generally not negative definite, Inequalities~\eqref{eq:thm_time_delay_LMI1} and~\eqref{eq:thm_time_delay_LMI2} are not convex in the state-feedback gain $K$ and hence, we cannot optimize over $K$ efficiently.
The following result provides a reformulation of Inequalities~\eqref{eq:thm_time_delay_LMI1}-\eqref{eq:thm_time_delay_LMI2} which is convex in $K$ for fixed variables $P_1^{-1}$, $R$ and hence, it can be used to search for a controller while simultaneously guaranteeing a bound on the corresponding MSI.

\begin{theorem}\label{thm:design}
If Assumption~\ref{ass:inertia} holds, then $\dot{x}(t)=Ax(t)+Bu(t)$ with controller~\eqref{eq:sampled_data_controller} is exponentially stable for any $\begin{bmatrix}A&B\end{bmatrix}\in\Sigma_{c}$ and any sampling sequence $\{t_k\}_{k=0}^{\infty}$ with $t_{k+1}-t_k\in(0,h]$, $t_0=0$, and $\lim_{k\to\infty}t_k=\infty$, if there exist ${\tilde{\lambda}_1>0}$, ${\tilde{\lambda}_2>0}$, ${Q_1\succ0}$, $R\succ0$, $Q_2$, and $Q_3$ with $Q_3+Q_3^\top\succ0$ such that $Q=\begin{bmatrix}Q_1&0\\Q_2&Q_3\end{bmatrix}$, $Q_R=\text{diag}(Q,R),\bar{Q}_R=\text{diag}(Q,R^{-1})$,~\eqref{eq:thm_time_delay_LMI3_dual}, and~\eqref{eq:thm_time_delay_LMI2_dual} hold.

\begin{figure*}
\vspace{2pt}
\begin{align}\label{eq:thm_time_delay_LMI3_dual}
\left[
\begin{array}{cc}
I&0\\
\begin{bmatrix}0&0&0\\I&-I&R\\0&0&-\frac{1}{2h}I\end{bmatrix}&
\begin{bmatrix}I&K^\top\\0&0\\0&0\end{bmatrix}\\\hline
0&I\\
\begin{bmatrix}0&I&0\end{bmatrix}&\begin{bmatrix}0&0\end{bmatrix}
\end{array}
\right]^\top
&\left[
\begin{array}{c|c}
\begin{matrix}
0&Q_R^\top\\Q_R&0
\end{matrix}&0\\\hline0&\tilde{\lambda}_1P_c
\end{array}
\right]
\left[
\begin{array}{cc}
I&0\\
\begin{bmatrix}0&0&0\\I&-I&R\\0&0&-\frac{1}{2h}I\end{bmatrix}&
\begin{bmatrix}I&K^\top\\0&0\\0&0\end{bmatrix}\\\hline
0&I\\
\begin{bmatrix}0&I&0\end{bmatrix}&\begin{bmatrix}0&0\end{bmatrix}
\end{array}
\right]\prec0\\\label{eq:thm_time_delay_LMI2_dual}
\left[
\begin{array}{cc}
I&0\\
\begin{bmatrix}0&0&0\\I&-I&0\\0&0&-\frac{1}{2}hI\end{bmatrix}&
\begin{bmatrix}I&K^\top\\0&0\\0&-hK^\top\end{bmatrix}\\\hline
0&I\\
\begin{bmatrix}0&I&0\end{bmatrix}&\begin{bmatrix}0&0\end{bmatrix}
\end{array}
\right]^\top
&\left[
\begin{array}{c|c}
\begin{matrix}
0&\bar{Q}_R^\top\\\bar{Q}_R&0
\end{matrix}&0\\\hline0&\tilde{\lambda}_2P_c
\end{array}
\right]
\left[
\begin{array}{cc}
I&0\\
\begin{bmatrix}0&0&0\\I&-I&0\\0&0&-\frac{1}{2}hI\end{bmatrix}&
\begin{bmatrix}I&K^\top\\0&0\\0&-hK^\top\end{bmatrix}\\\hline
0&I\\
\begin{bmatrix}0&I&0\end{bmatrix}&\begin{bmatrix}0&0\end{bmatrix}
\end{array}
\right]\prec0
\end{align}
\noindent\makebox[\linewidth]{\rule{\textwidth}{0.4pt}}
\end{figure*}
\end{theorem}
\textbf{Proof}$\quad$
Applying the dualization lemma~\cite[Lemma 4.9]{scherer2000linear} and performing straightforward algebraic manipulations, we infer that~\eqref{eq:thm_time_delay_LMI3_dual} is equivalent to
\begin{align}\label{eq:thm_time_delay_LMI3}
\left[
\begin{array}{cc}
I&0\\
\begin{bmatrix}0&I&0\\0&-I&0\\0&R&-\frac{1}{2h}I\end{bmatrix}&\begin{bmatrix}0\\I\\0\end{bmatrix}\\\hline
0&I\\
\begin{bmatrix}I&0&0\\K&0&0\end{bmatrix}&\begin{bmatrix}0\\0\end{bmatrix}
\end{array}
\right]^\top
\left[
\begin{array}{c|c}
\begin{matrix}
0&P_R^\top\\P_R&0
\end{matrix}&0\\\hline0&\lambda_1\tilde{P}_c
\end{array}
\right]
\left[
\begin{array}{cc}
\star&\star\\\star&\star\\\hline
\star&\star\\\star&\star
\end{array}
\right]\prec0
\end{align}
with $\lambda_1=\frac{1}{\tilde{\lambda}_1},P=\begin{bmatrix}P_1&0\\P_2&P_3\end{bmatrix}\coloneqq Q^{-1},P_R\coloneqq\text{diag}(P,R^{-1})$, where invertibility of $Q$ follows from $Q_1\succ0,Q_3+Q_3^\top\succ0$.
Similar to the proof of Theorem~\ref{thm:analysis}, the full-block S-procedure~(\cite{scherer2001lpv}) together with~\eqref{eq:AB_bound_primal} implies that
\begin{align*}
&\begin{bmatrix}P_1&0&0\\P_2&P_3&0\\0&0&R^{-1}\end{bmatrix}^\top
\begin{bmatrix}0&I&0\\A+BK&-I&0\\0&R&-\frac{1}{2h}I\end{bmatrix}\\
&+\begin{bmatrix}0&I&0\\A+BK&-I&0\\0&R&-\frac{1}{2h}I\end{bmatrix}^\top
\begin{bmatrix}P_1&0&0\\P_2&P_3&0\\0&0&R^{-1}\end{bmatrix}\\
=&\begin{bmatrix}P_2^\top(A+BK)+(A+BK)^\top P_2&\star&\star\\
P_1-P_2+P_3^\top(A+BK)&-P_3-P_3^\top&\star\\
0&I&-\frac{1}{h}R^{-1}\end{bmatrix}\\
\prec&\>0
\end{align*}
for any $\begin{bmatrix}A&B\end{bmatrix}\in\Sigma_c$.
Applying the Schur complement to the right-lower block, this in turn implies~\eqref{eq:time_delay1}.
Similarly, using the dualization lemma~\cite[Lemma 4.9]{scherer2000linear},~\eqref{eq:thm_time_delay_LMI2_dual} implies~\eqref{eq:thm_time_delay_LMI2} which in turn implies~\eqref{eq:time_delay2} by Theorem~\ref{thm:analysis}, thus concluding the proof.
$\hfill\blacksquare$

It is difficult to render Inequalities~\eqref{eq:thm_time_delay_LMI3_dual}--\eqref{eq:thm_time_delay_LMI2_dual} convex in all decision variables since i) they contain both $R$ and $R^{-1}$ and ii) the state-feedback gain $K$ is multiplied by both $Q_1$ and $R^{-1}$.
Nevertheless, for fixed $Q_1\succ0$ and $R\succ0$, they are convex in the remaining decision variables and thus, Theorem~\ref{thm:design} allows us to use measured data for designing controllers which lead to a possibly large MSI bound $h$.
While fixing the matrices $Q_1,R$ leads to conservatism, a practical remedy is to alternate between solving the LMIs in Theorems~\ref{thm:analysis} and~\ref{thm:design} for fixed $K$ and $Q_1,R$, respectively.
The proof of Theorem~\ref{thm:design} exploits dualization arguments to render the involved matrix inequalities convex in $K$, similar to the (discrete-time) results in~\cite{berberich2020combining}.
In particular, the dualization lemma directly implies that~\eqref{eq:thm_time_delay_LMI2} and~\eqref{eq:thm_time_delay_LMI2_dual} are equivalent.
In order to render Inequality~\eqref{eq:thm_time_delay_LMI1} convex in $K$, we additionally need to perform a Schur complement w.r.t. $R$.

\begin{remark}
In addition to being purely data-driven methods, Theorems~\ref{thm:analysis} and~\ref{thm:design} also extend existing model-based methods for stability analysis and controller design for uncertain sampled-data systems.
Many existing time-delay approaches such as~\cite{fridman2010refined,gao2010robust,seuret2012novel} can only handle polytopic uncertainty and verify the nominal LMIs on all vertices of the uncertainty set, which can become intractable for too many vertices (e.g., in high-dimensional systems).
Alternatively, the approach in~\cite{orihuela2010delay} allows for simple quadratic uncertainty bounds but only stability analysis conditions are provided and no controller design is considered.
On the other hand, the presented approach utilizes powerful robust control tools to employ the data-dependent uncertainty bound derived in Theorem~\ref{thm:cont_param}.
Additionally, it is straightforward to include prior model knowledge in order to reduce conservatism and improve the MSI bound:
Suppose the LTI system~\eqref{eq:sys_cont} is given as a linear fractional transformation (compare~\cite{zhou1996robust}) 
\begin{align}\label{eq:LFT}
\dot{x}(t)&=Ax(t)+Bu(t)+B_ww(t),\\\nonumber
z(t)&=Cx(t)+Du(t),\\\nonumber
w(t)&=\Delta z(t),
\end{align}
where all matrices except for the uncertainty $\Delta$ are known.
Following the arguments in~\cite{berberich2020combining}, all results in this paper can be directly extended to compute MSI bounds and design controllers for~\eqref{eq:LFT} robustly for all uncertainties $\Delta$ which are consistent with the measured data (potentially affected by noise) and the prior knowledge given by the matrices $A,B,C,D,B_w$.
This also allows us to include prior knowledge on uncertainty bounds or structure such as
\begin{align*}
\Delta=\mathrm{diag}(\Delta_1,\dots,\Delta_{\ell}),\quad
\begin{bmatrix}\Delta_i^\top\\I\end{bmatrix}^\top
P_i
\begin{bmatrix}\Delta_i^\top\\I\end{bmatrix}\succeq0
\end{align*}
for some matrices $P_i$, which can lead to a dramatic reduction of uncertainty and therefore to an increase of the computed MSI bound.
Finally, for the special case that only prior knowledge and no data is available, this leads to a model-based robust time-delay approach for uncertain systems with an improved flexibility if compared to the existing literature mentioned above. 
\end{remark}

\begin{remark}
As an obvious alternative to the results in this paper, one could estimate the matrices $A_{\tr}$, $B_{\tr}$ together with uncertainty bounds and then apply the model-based time-delay results for uncertain systems in~\cite{fridman2010refined,gao2010robust,seuret2012novel,orihuela2010delay}.
However, as discussed above, these results mainly focus on polytopic uncertainty descriptions.
It is possible to construct a polytope which is guaranteed to contain the unknown system matrices $A_{\tr}$, $B_{\tr}$ using measured data affected by noise, e.g., via set membership estimation (\cite{milanese1991optimal}).
However, the number of vertices of this polytope grows rapidly with the system dimension and the number of data points and hence, if the polytope is not over-approximated (leading to conservatism), then this approach quickly becomes intractable. 
\end{remark}

Finally, to the best of our knowledge, the presented work provides the first results to compute MSI bounds and design controllers with a possibly high MSI based only on measured data.
A highly interesting issue for future research is whether and, if so, how further approaches to sampled-data control under aperiodic sampling as, e.g., surveyed in~\cite{hetel2017recent} can be translated into a purely data-driven formulation.
To this end, following the input-output approach (compare~\cite{mirkin2007some}), which relies on an integral quadratic constraint description (\cite{megretski1997system}) of the delay operator, in such a scenario is immediate, simply by considering both the unknown system matrices and the delay as uncertainties in a robust control formulation (compare~\cite{berberich2020combining} for how additional uncertainty bounds can be combined with a data-dependent bound of the form~\eqref{eq:thm_cont_param} to design robust controllers).

\section{Example}\label{sec:example}

In the following, we illustrate the practical applicability of the proposed approach by applying it to a standard example from~\cite{zhang2001stability}.
To be precise, we consider
\begin{align}\label{eq:sys_ex}
A_{\tr}=\begin{bmatrix}0&1\\0&-0.1\end{bmatrix},\quad
B_{\tr}=\begin{bmatrix}0\\0.1\end{bmatrix},\quad B_d=I,
\end{align}
where $A_{\tr},B_{\tr}$ are \emph{unknown}.
We assume that $N=100$ measurements $\{\dot{x}(\tau_k),x(\tau_k),u(\tau_k)\}_{k=1}^N$ are available with (not equidistant) sampling intervals
\begin{align*}
&\tau_{k+1}-\tau_k=1.5\quad\text{if}\>k\in[1,49],\\
&\tau_{k+1}-\tau_k=3\quad\text{if}\>k\in[50,99],
\end{align*}
where the input generating the data is sampled uniformly from $u(\tau_k)\in[-1,1]$.
These measurements are perturbed by a disturbance $\{\hat{d}(\tau_k)\}_{k=1}^N$ sampled uniformly from $\hat{d}(\tau_k)^\top\hat{d}(\tau_k)\leq\bar{d}^2$, $k=1,\dots,N$, for some $\bar{d}>0$, and we assume the (valid) noise bound~\eqref{eq:noise_bound} with $Q_d=-I$, $S_d=0$, $R_d=\bar{d}^2NI$.
We use MATLAB together with YALMIP~(\cite{lofberg2004yalmip}) and MOSEK~(\cite{MOSEK15}) for computations.
We first apply Theorem~\ref{thm:analysis} to compute an MSI bound for a given controller $K=-\begin{bmatrix}3.75&11.5\end{bmatrix}$ which is also considered in~\cite{fridman2010refined} and for which solving the model-based LMIs~\eqref{eq:time_delay1}-\eqref{eq:time_delay2} leads to a maximum $h$ of $1.62$.
On the other hand, the MSI bounds computed based on the available data without any model knowledge via Theorem~\ref{thm:analysis} are listed in Table~1 for varying noise levels.
As expected, the MSI bound which can be guaranteed decreases for increasing noise levels since the corresponding set $\Sigma_c$ of systems consistent with the data grows and hence, the sampling frequency required for stability increases.

Since the bound~\eqref{eq:noise_bound} with $Q_d=-I$, $S_d=0$, $R_d=\bar{d}^2NI$ is not a tight representation of $\hat{d}(\tau_k)^\top\hat{d}(\tau_k)\leq\bar{d}^	2$, the size of $\Sigma_c$ need not decrease in general when the data length $N$ increases.
This phenomenon has been pointed out and analyzed in more detail in the recent literature, compare, e.g.,~\cite{martin2020dissipativity,berberich2020combining,bisoffi2021trade}.
Nevertheless, for the application of Theorems~\ref{thm:analysis} and~\ref{thm:design} in the above example, we observe that increasing the data length $N$ has (approximately) a similar effect as reducing the noise bound $\bar{d}$, i.e., a longer data trajectory allows us to compute larger MSI bounds in most cases.

\begin{table}\label{tab:ex_analysis}
\begin{center}
\begin{tabular}{c|ccccccc}
noise bound $\bar{d}$&$0.001$&$0.005$&$0.01$&$0.02$&$0.03$&$0.04$&$0.05$\\\hline
MSI bound $h$&
$1.59$&$1.49$&$1.38$&$1.17$&$1$&$0.86$&$0.67$
\end{tabular}
\vskip6.5pt
\caption{MSI bound $h$ for~\eqref{eq:sys_ex} with 
$K=-\big[3.75\quad11.5\big]$, computed via Theorem~\ref{thm:analysis}.}
\end{center}
\end{table}

Further, we optimize for a controller gain $K$ with a possibly large MSI bound by iteratively optimizing over the LMIs in Theorems~\ref{thm:analysis} and~\ref{thm:design} (fixing $K$ and $Q_1,R$, respectively) while gradually increasing $h$.
We initialize this iterative approach via the controller $K=-\begin{bmatrix}3.75&11.5\end{bmatrix}$.
In Table~2, the MSI bounds $h$ are displayed which hold for the corresponding state-feedback gains that emerged from this iteration.
Note that, by choosing a different controller, the MSI bound can be increased significantly if compared to the controller given above.

\begin{table}\label{tab:ex_design}
\begin{center}
\begin{tabular}{c|ccccccc}
noise bound $\bar{d}$&$0.001$&$0.005$&$0.01$&$0.02$&$0.03$&$0.04$&$0.05$\\\hline
MSI bound $h$&
$142.6$&$28.5$&$13.8$&$6.3$&$4$&$2.9$&$2.2$
\end{tabular}
\vskip6.5pt
\caption{MSI bound $h$ for~\eqref{eq:sys_ex} with $K$ computed via an iteration using Theorems~\ref{thm:analysis} and~\ref{thm:design}.}
\end{center}
\end{table}

\bibliography{Literature}  

\section*{Appendix A: Error bounds for state derivative estimation}
In the following, we show how data $\{x(\tau_k),u(\tau_k)\}_{k=1}^N$ can be used to derive an upper bound on the derivative estimation error resulting from a simple Euler discretization.
For simplicity, we assume that the available data are sampled equidistantly, i.e., $\tau_{k+1}-\tau_k=h>0$ for $k=1,\dots,N-1$.
Moreover, we assume that the data are not affected by noise, i.e., $\hat{d}(\tau_k)=0$ for $k=1,\dots,N$, but we note that the case of noisy data can be handled analogously using the bound~\eqref{eq:noise_bound}.
Finally, we assume that norm bounds on the unknown system matrices are available, i.e., $\lVert A_{\tr}\rVert_2\leq\bar{a}$ and $\lVert B_{\tr}\rVert_2\leq\bar{b}$ with known $\bar{a},\bar{b}>0$.

For $k=1,\dots,N-1$, we construct an estimate of $\dot{x}(\tau_k)$ via
\begin{align}\label{eq:app_A_0}
\underline{\dot{x}}(\tau_k)=\frac{x(\tau_{k+1})-x(\tau_k)}{h}.
\end{align}
Using the explicit solution formula for~\eqref{eq:sys_cont}, we have
\begin{align}\label{eq:app_A_1}
x(\tau_{k+1})=e^{A_{\tr}h}x(\tau_k)+\int_0^he^{A_{\tr}s}\mathrm{d}sB_{\tr}u(\tau_k).
\end{align}
Taylor's Theorem implies
\begin{align}\label{eq:app_A_2}
e^{A_{\tr}h}&=I+A_{\tr}h+\frac{1}{2}(A_{\tr}q_1)^2,\\\label{eq:app_A_3}
\int_0^he^{A_{\tr}s}\mathrm{d}s&=hI+\frac{1}{2}A_{\tr}h^2+\frac{1}{2}A_{\tr}^2\int_0^hq_2(s)^2\mathrm{d}s
\end{align}
for some $q_1\in[0,h]$, $q_2(s)\in[0,s]$.
Plugging~\eqref{eq:app_A_1}--\eqref{eq:app_A_3} into~\eqref{eq:app_A_0}, we obtain
\begin{align*}
\underline{\dot{x}}(\tau_k)=&\dot{x}(\tau_k)+\frac{1}{2h}(A_{\tr}q_1)^2x(\tau_k)\\
&+\left(\frac{1}{2}A_{\tr}h+\frac{1}{2h}A_{\tr}^2\int_0^hq_2(s)^2\mathrm{d}s\right)B_{\tr}u(\tau_k).
\end{align*}
Using $q_2(s)\leq s$, this implies the error bound $\lVert\underline{\dot{x}}(\tau_k)-\dot{x}(\tau_k)\rVert_2\leq\bar{x}$ with
\begin{align*}
\bar{x}\coloneqq \frac{\bar{a}h}{2}\left(\bar{a}\lVert x(\tau_k)\rVert_2+\left(1+\frac{\bar{a}h}{3}\right)\bar{b}\lVert u(\tau_k)\rVert_2\right).
\end{align*}

\section*{Appendix B: Satisfaction of Assumption~\ref{ass:inertia}}        
In the following, we show that Assumption~\ref{ass:inertia} holds under the following conditions: $S_d=0$, $Z$ has full row rank, $B_d$ is invertible, and the disturbance $\hat{D}$ generating the data satisfies $\hat{D}Q_d\hat{D}^\top+R_d\succ0$.
If $Z$ has full row rank, then $Q_c=ZQ_dZ^\top\prec0$ such that the following can be easily derived (compare~\cite{scherer2000linear})
\begin{align*}
\begin{bmatrix}I&0\\-S_c^\top Q_c^{-1}&I\end{bmatrix}P_c
\begin{bmatrix}I&0\\-S_c^\top Q_c^{-1}&I\end{bmatrix}^\top=
\begin{bmatrix}Q_c&0\\0&R_c-S_c^\top Q_c^{-1}S_c
\end{bmatrix}.
\end{align*}
Thus, Assumption~\ref{ass:inertia} holds if $Z$ has full row rank and $R_c\succ S_c^\top Q_c^{-1} S_c$.
Inserting $S_d=0$ and using the definition of $Q_c$, $S_c$, $R_c$ (compare~\eqref{eq:bound_def}), the latter inequality is equivalent to
\begin{align}\label{eq:app_1}
B_dR_dB_d^\top\succ\dot{X}\left(Q_dZ^\top(ZQ_dZ^\top)^{-1}ZQ_d-Q_d\right)\dot{X}^\top.
\end{align}
By assumption, the data contained in the matrices $\dot{X}$ and $Z$ are generated by the disturbance $\hat{D}$, i.e., $\dot{X}=\begin{bmatrix}A_{\tr}&B_{\tr}\end{bmatrix}Z+B_d\hat{D}$.
Plugging this into~\eqref{eq:app_1} and using that
\begin{align*}
\begin{bmatrix}A_{\tr}&B_{\tr}\end{bmatrix}Z\left(Q_dZ^\top(ZQ_dZ^\top)^{-1}ZQ_d-Q_d\right)=0,
\end{align*}
Inequality~\eqref{eq:app_1} is equivalent to
\begin{align}\label{eq:app_2}
B_dR_dB_d^\top\succ B_d\hat{D}\left(Q_dZ^\top(ZQ_dZ^\top)^{-1}ZQ_d-Q_d\right)\hat{D}^\top B_d^\top.
\end{align}
Since $B_d$ is invertible, the condition $\hat{D}Q_d\hat{D}^\top+R_d\succ0$ implies
\begin{align}\label{eq:app_3}
B_dR_dB_d^\top\succ-B_d\hat{D}Q_d\hat{D}^\top B_d^\top.
\end{align}
Using~\eqref{eq:app_3}, Inequality~\eqref{eq:app_2} holds if
\begin{align}\label{eq:app_4}
B_d\hat{D}Q_dZ^\top(ZQ_dZ^\top)^{-1}ZQ_d \hat{D}^\top B_d^\top\preceq0,
\end{align}
which is true since $(ZQ_dZ^\top)^{-1}\prec0$.
To conclude, we have shown that, under the proposed conditions (i)--(iii), Inequality~\eqref{eq:app_1} and hence Assumption~\ref{ass:inertia} holds.
\end{document}